\documentclass{article}
\usepackage[utf8]{inputenc}
\usepackage{todonotes}
\usepackage{authblk}

\title{State-space reduction techniques exploiting specific constraints for quantum search \\ Application to a specific job scheduling problem}
\author[1]{Rodolphe Griset}
\author[2]{Ioannis Lavdas}
\author[3]{Jiri Guth Jarkovsky}
\affil[1]{EDF R\&D}
\affil[2]{Arnold-Sommerfeld-Center for Theoretical Physics\\
Ludwig-Maximilians-Universit\"at, 80333 M\"unchen, Germany\footnote{As of September 2023, Welinq, Paris}}
\affil[3]{IQM Germany GmbH, Georg-Brauchle-Ring 23-25, 80992 M\"unchen, Germany}

\date{\today}

\usepackage[numbers]{natbib}
\usepackage{graphicx}
\usepackage{tikz}
\usepackage{todonotes}
\usepackage{amsmath}
\usepackage{subcaption}
\usepackage{dsfont}
\usepackage{tcolorbox}
\usepackage{verbatim} 
\usepackage{amssymb} 
\usepackage{pgfplots}
\usepackage{tcolorbox}

\usepackage{url} 

\definecolor{lightgray}{gray}{0.75}
\usetikzlibrary{quantikz}
\newtheorem{assumption}{Assumption}
\begin{document}

\maketitle

\textbf{Abstract.} Quantum search has emerged as one of the most promising fields in quantum computing. State-of-the-art quantum search algorithms enable the search for specific elements in a distribution by monotonically increasing the density of these elements until reaching a high dentisity. This kind of algorithms demonstrate a theoretical quadratic speed-up on the number of queries compared to classical search algorithms in unstructured spaces. Unfortunately, the major part of the existing literature applies quantum search to problems which size grows exponnentialy with the input size without exploiting any specific problem structure, rendering this kind of approach not exploitable in real industrial problems. In contrast, this work proposes exploiting specific constraints of scheduling problems to build an initial superposition of states with size almost quadraticaly increasing as a function of the problem size. This state space reduction, inspired by the quantum walk algorithm, constructs a state superposition corresponding to all paths in a state-graph embedding spacing constraints between jobs.
Our numerical results on quantum emulators highlights the potential of state space reduction approach, which could lead to more efficient quantum search processes by focusing on a smaller, more relevant, solution space.

\section{Introduction}
Since the introduction of the renowned Grover search algorithm (\cite{grover1996fast}), which demonstrated a quadratic speedup in searching for a solution in an unsorted database, quantum search has emerged as one of the most promising fields of quantum computing. This family of algorithms relies on a so-called oracle function to identify elements of interest in a space of possible solutions. The algorithm begins with an initial state superposition, classically a uniform one which is easy to build. It then iteratively increases the density of marked elements by switching their phase and using this switch to rotate from the initial state superposition to a superposition where the marked elements have high density. The initial Grover algorithm necessitates a specific number of iterations to circumvent the soufflé problem, corresponding to stop the algorithm before reaching the desired distribution or at the opposite overshoot it and beginning to reduce the density of marked element in the distribution \cite{brassard1997searching}. Over time, algorithms based on Grover have been extended to search for several marked elements (\cite{brassard2002quantum}) and to avoid the soufflé problem by adapting the rotation to only increase the density of marked elements at each iteration (\cite{yoder2014fixed}). While Grover-based algorithms are conceptually appealing, there is little work applying these algorithms to real-world industrial use cases (\cite{botsinis2018quantum},\cite{zajac2022towards}). In fact, this type of algorithm is primarily employed to identify the optimal or near-optimal solution to optimization problems, which necessitate a substantial depth and a high qubit fidelity not available in current NISQ (Noisy Intermediate Scale Quantum) computers. Secondly, the majority of operational research problems solved in the industry have a state-space of an exponential size with the input data, for which a quadratic speed-up might be insufficient to justify moving from a classical high performance computer (HPC) to a quantum computer (QC) in the near future.

In this article, we propose to exploit the fact that quantum search algorithms do not impose specific conditions on the initial state superposition. By applying these algorithms to a state superposition that already ensures some constraints of the problem are satisfied, we can reduce the complexity of the quantum search. Constraints make some possibilities in the initial space infeasible, thus the number of elements with a non-zero density in this state superposition, which we will call the \textit{reduced space state}, is lower than the total number of possibilities in the solution space. This reduction decreases the quantum search complexity. An algorithm exploiting this principle divides the problem’s constraints into two categories: the first are the one which have an exploitable structure which can be handle directly through the construction of the reduced space state, while oracle marks elements of the second category, without exploitable structure, to be targeted by the quantum search approach.

We introduce such categorization for a specific kind of scheduling problem involving a set of machines and a set of jobs to be scheduled on each machine. This problem exhibits specific time window structures for jobs and spacing constraints between jobs on the same machines. Jobs from different machines share resources, which limits the number of jobs that can be performed in parallel. In the classical world, resource constraints render the problem difficult to solve in general, while spacing constraints may be efficiently embedded in an extended state graph formulation \cite{pira2020novel}. Our approach exploits this extended state graph to build the reduced space state before applying quantum search to handle resource constraints.

We developped a quantum walk-inspired method to build the reduced space state. To the best of our knowledge, this approach of state-space reduction via quantum walk for quantum search is new in the literature, although similar ideas have been studied to improve the initial state of variational optimization algorithms. In particular, the textbook \emph{quantum approximate optimization algorithm} (QAOA) \cite{farhi2014quantum} also starts with the superposition of all computational basis states and, through the application of specific unitaries, tries to find the states that best solve an optimization problem. Hadfield et. al \cite{hadfield2019from} proposed an evolution of QAOA, the \emph{quantum alternating operator ansatz} (QAOAnsatz), which works similarly but starts with the superposition of only the states satisfying a fixed constraint\footnote{Common constraints include, for example, a fixed Hamming weight (number of qubits in the $\ket{1}$ state) or a one-hot encoding (groups of qubits of which exactly one is in the $\ket{1}$ state.)}. Correspondingly, the QAOAnsatz algorithm is designed to maintain states holding this constraint throughout the computation. This kind of technic has been also used to biase the initial state superposition toward good solutions using classical heuristics \cite{van2021quantum} which act as a kind of warm-start for variationnal algorithms.




Our contributions are the following :
\begin{itemize}
\item We design a quantum search approach for a specific scheduling problems and show that the quadratic speed-up promised by this kind of approaches will not be sufficient compared to the exponential search space size. 
\item We propose a quantum walk-inspired scheme to reduce the state-space from a quantum superposition of states defined by time windows and spacing constraints between jobs, before applying the quantum search algorithm to handle resource constraints.
\item As a proof of concept, we implemented this scheme on a simplified use case using a quantum emulator and compared the impact of problem size on both full and reduced search approaches, showcasing the high potential of state-space reduction.
\end{itemize}
The remaining part of this article is organized as follows : section \ref{sec:IndustrialUseCase} presents the industrial motivations and the simplified problem, \ref{sec:QSearch} presents the application of a quantum search approach to this simplified use case. Section \ref{sec:SSR} introduces the proposed state-space reduction technique inspired from discrete quantum walk approach. Finally, section \ref{sec:Results} gives numerical results demonstrating the importance of the state-space reduction technique for quantum search problems.

\section{Problem description} \label{sec:IndustrialUseCase}

\subsection{Generic use case}

We consider a specific variant of job scheduling satisfability problem considering a set of $I$ machines which have to perform $K$ jobs each, over a set of time steps $T$. Time windows are associated with each job, following a specific structure : 
\begin{itemize}
\item A time window $[E_{i},D_{i}]$, is associate with the first job from machine $i \in I$. We call the \textit{Offset} $O_i$ of the machine $i$, the difference between the earliest dates of the machine and the overall earliest time of all first job on all machines, namely $O_i = E_i - \min_i E_i$.
\item A minimum time ($\underline{T}_i$) and a maximum time ($\overline{T}_i$) are required between two consecutive jobs of the same machine. Those minimum and maximum times combined with the time window of the first job of the machine define the time window structure of further jobs, namely $[E_{ik},D_{ik}]=[E_i + k.\underline{T}_i,D_i + k.\overline{T}_i]$.
\end{itemize}

Note that, the time window of jobs corresponding to index $k$ is of size $C.k+1$. Finaly, a set of resources constraints $c \in C$ limit the number of jobs allowed to be perform in parallel on different machines. In general, an objective function is associated with this kind of problem such as total completion time or minimum tardiness, this kind of optimisation problems are NP-hard in the strong sense for more than one machine \cite{blazewicz1983scheduling}.   

\subsection{Simplified case}

Given emulator and current machine limitations, we consider only satisfability problems where we aim at finding a set of feasible solutions. The motivation for considering this problem will be detailed in section \ref{subsec:InduMotiv}. Moreover, we introduce some simplifying assumptions in order to reduce circuit complexity :

\begin{assumption} \label{ass:Spacing}
\textbf{Spacing constraints} : Spacing constraints correspond to the following expression: 
$$  0 \leq X_{i,k+1} - X_{i,k} < C$$ 
\end{assumption}

\begin{figure}[!ht]
\centering
\resizebox{.8\linewidth}{!}{
\begin{tikzpicture}[shorten >=1pt,->]
  \tikzstyle{vertex}=[circle,fill=black!25,minimum size=17pt,inner sep=0pt]
  
  \node[vertex](G-s) at (0,0) {s};
  
  \foreach \y in {0, 1, 2, 3}
    		\node[vertex] (G-1-\y) at (5,2*\the\numexpr\y-1) {$\y$};
  
  \foreach \y in {0,1,2,3,4,5,6}
  	\node[vertex] (G-2-\y) at (10,2*\the\numexpr\y-2) {$\y$};

  \draw[red, dashed, thick] ($(G-1-3.north west)+(-0.2,0.2)$) rectangle ($(G-1-0.south east)+(0.2,-0.2)$);
  \node[red, anchor=east] at ($(G-1-3.north west)+(-0.2,0)$) {Dates job 1};
  
  \draw[blue, dashed, thick] ($(G-2-6.north west)+(-0.2,0.2)$) rectangle ($(G-2-0.south east)+(0.2,-0.2)$);
  \node[blue, anchor=east] at ($(G-2-4.north east)+(3,0)$) {Dates job 2};

  	\foreach \node in {0,1,2,3}
  	{
    		\draw (G-s) -- (G-1-\node);
    		\foreach \othernode in {\the\numexpr\node+3,\the\numexpr\node+2,\the\numexpr\node+1, \the\numexpr\node}
    		{
    			\draw (G-1-\node) -- (G-2-\othernode);
    		}	
    	}

\end{tikzpicture}}
\caption{Graphic representation of the simplified instances structure and label}\label{fig:GraphExemple}
\end{figure}
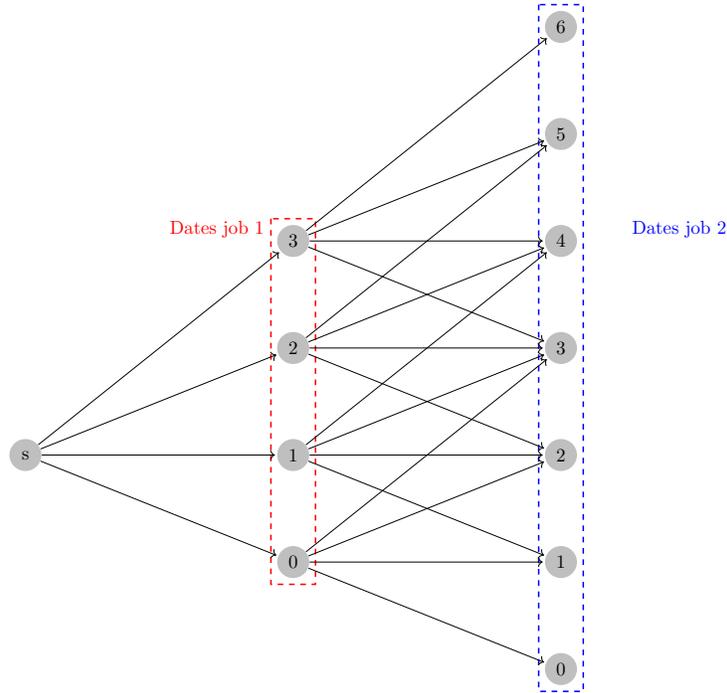

Figure \ref{fig:GraphExemple} presents a graphical example of the structure for $I = K = 2$ and $C=4$. Each
 node represents one possible date and each edge corresponds to possible combinations of dates allowed by the spacing constraints.  This example will be use to detail our algorithms in the following of this article

\begin{assumption}\label{ass:ResourceUnitary}
    \textbf{Unitary capacity :} We assume the capacity of the resource constraint is one which means that the resource constraint do not allow to choose the same dates for two jobs of different machines.
\end{assumption}

\begin{assumption} \label{ass:ResourceCT}
\textbf{Resource constraints by job index} : We assume that resource constraints only impact jobs with the same index $k$. 
\end{assumption}

Note that these assumptions are not strictly necessary and can be relaxed by using a higher number of qubits and gates. They are practical limitations used to simplify the presentation of our approaches and to enable simulations, but there is no theoretical bottleneck here. Assumption \ref{ass:Spacing} simplifies the structure of time windows and spacing constraints, allowing them to be represented on a regular graph, which facilitates the design of a quantum circuit that exploits this structure. Assumption \ref{ass:ResourceUnitary} allows resource constraints to be checked with a simple CNOT operation between qubits, but it can be easily extended to greater values by using additional qubits to store information about the number of jobs in progress and compare it to a threshold. Finally, Assumption \ref{ass:ResourceCT} facilitates operations by only comparing qubits representing the same job index, but it can be extended to compare the dates of all jobs.

\subsection{Industrial motivation} \label{subsec:InduMotiv}

This problem is a refined model of the outage planning problem \cite{griset2021nuclear}, which additionally associates a non-convex cost function with outage dates. The outage planning problem considers a set of machines that must be stopped periodically and aims to find a set of outage dates for each unit under resource constraints. These constraints limit the number of outages from different machines that can be performed simultaneously and impose spacing constraints between outages of the same unit, modeling specific requirements such as fuel management or periodic maintenance.
Maintenance planning problems are crucial for energy producers as they shape the entire electricity market for short-term issues. Generally, these are complex, large-scale problems involving scheduling constraints on outages, various technical constraints related to the machines, and a stochastic component. The primary objective is often to meet demand at the lowest possible cost in an uncertain future. Despite extensive research on these problems in the literature \cite{kobbacy2008maintenance}, they remain difficult to be solved optimally and are mainly addressed using meta-heuristic approaches in the industry. This is where quantum computing might become a game changer by generating a set of promising solutions that can be used as initial information for classical algorithms or experts in the field. Indeed, it’s not always possible to give a formal description of all the constraints and the cost function associated with these kinds of problems. In this context, our idea is to exploit quantum computing not to solve a problem to optimality but to generate a pool of candidate solutions to be used as input for a more precise classical algorithm or to be analyzed by an expert. The idea of this work is to build a first step in this direction by designing a quantum search approach to generate feasible solutions to these kinds of problems. Precisely for that we consider in this paper the feasibility problem, aiming to find feasible solutions to the maintenance planning problem.

\section{Quantum search approach}\label{sec:QSearch}

\subsection{State of the art of quantum search algorithms} \label{sec:SOTA}

The first quantum search algorithm introduced was Grover’s algorithm \cite{grover1996fast}. This algorithm relies on a circuit, called oracle, that “marks” the target state by flipping its quantum phase. The oracle is used in conjunction with the \emph{Grover diffusion operator}, which reflects all states about the average amplitude of the superposition. The algorithm iteratively applies the oracle, which makes the average amplitude of the superposition lower than the phase of the non-marked elements but greater than the (negative) phase of the marked element. The diffusion operator then reflects all states about this average, thereby amplifying the amplitude of the target state within the superposition. After approximatively, in $\mathcal{O}\left( \sqrt{N}\right)$ iterations, the target state becomes the dominant state of the superposition, making it highly probable to be measured. In contrast, the best classical algorithm for unstructured search requires $\mathcal{O}(N)$ queries to a classical oracle. Under standard complexity-theoretic assumptions, Grover’s algorithm is considered optimal for unstructured search jobs \cite{bennett1997strengths}. Despite its advantages, Grover’s algorithm has several limitations. These have been addressed by more advanced search algorithms that extend Grover’s approach in various ways. One such algorithm is \emph{amplitude amplification} \cite{brassard2002quantum}. In amplitude amplification, the search space is divided into two orthogonal subspaces: a “bad” subspace and a “good” subspace corresponding to $M$ solutions. Assuming an oracle can mark states in the good subspace, the algorithm iteratively applies the oracle and the diffusion operator, similar to Grover’s algorithm. This process amplifies the amplitude of the good states while decreasing the amplitude of the bad states. After approximately $\mathcal{O}\left( \sqrt{ N / M }\right)$ iterations, the superposition predominantly consists of good states, ensuring that a measurement will yield one of them.

If the number of good states is not known beforehand, amplitude amplification in its simple form cannot be used. Without knowing $M$, we cannot determine the number of Grover iterations needed to amplify the good states. Iterating too few times will leave many bad states in the superposition, while iterating too many times will start having the opposite effect, amplifying the bad states and suppressing the good states. This issue is addressed by \emph{fixed-point search} algorithms \cite{yoder2014fixed}, which are designed to converge towards the good subspace. Fixed-point search algorithms are based on the \emph{quantum singular value transformation} (QSVT) \cite{martyn2021grand}. Theoretically, QSVT-based searches provide the same quadratic speed-up as Grover’s algorithm. In practice, the speed-up depends on the user's uncertainty about the number of marked elements. Assuming at least $M$ marked elements, the algorithm can be run with $\mathcal{O}\left( \sqrt{N/M} \right)$ steps to find a marked element with high probability (if the assumption was correct). Classical brute-force search is expected to find a marked element in $N/m$ steps, where $m$ is the true (unknown) number of marked elements, so the true speed-up depends on the relation of $M$ and $m$.

In the extreme case when most of the states are good, amplitude amplification can be used for deleting the bad states with super-exponential speedup \cite{liu2013deleting}. Classically, this job would take $\mathcal{O}(N)$ steps, but with amplitude amplification this can be performed instead in $\mathcal{O}(1)$ steps. 


\subsection{Direct application of quantum search}

\begin{figure}[ht!]
    \begin{quantikz}
    \lstick{Qubits \\ machine 1}  
     & \gate{H} &\gate{\begin{array}{c}Feasible Path \\ Oracle \end{array}} & \gate[wires=2]{\begin{array}{c}Resource \\ Constraints \\\ Oracle \end{array}} & \gate[wires=2]{\begin{array}{c}Amplitude \\ amplification \end{array}}
    & \rstick{}\qw \\
    \lstick{Qubits \\ machine 2}  & \gate{H} &
    \gate{\begin{array}{c}Feasible Path \\ Oracle \end{array}}  & &
    & \rstick{} \qw \\
    \end{quantikz}
    \caption{Full quantum search algorithm principle}
    \label{fig:FS_schema}
\end{figure}
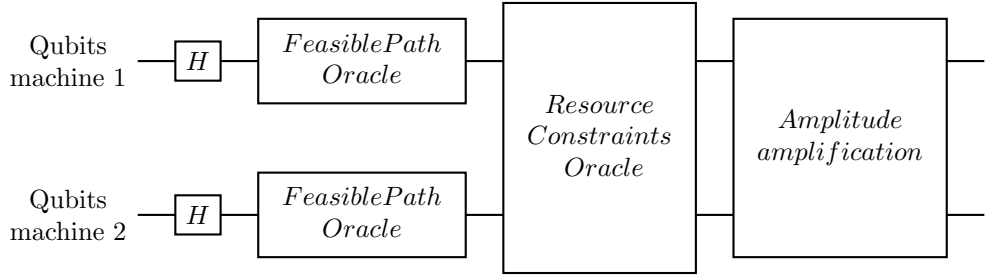

This section introduces the implementation of a quantum search for the simplified use case. To use this algorithm, we need to define an encoding that maps the solutions of the problem to qubit combinations and oracles, which will mark feasible solutions to apply amplitude amplification. Figure \ref{fig:FS_schema} illustrates the principle of the algorithm, which starts by putting the system in a superposition of all possible qubit combinations. Then, oracles are applied to mark the correct elements before applying the amplitude amplification algorithm.

\subsubsection{Encoding} \label{subsubsec:Encoding}
Resource constraints require comparing the absolute job positions of different machines. Therefore, we need to account for the offset between machines in the chosen encoding. Let $O$ be the maximum offset between machines, i.e $O = \max_i(D_{i0}) - \min_{i}(E_{i0}) $. We need $log_2(O)$ qubits to store the offset of each machine. Then, for each job index $k$, we associate a binary sequence which gives the date of the corresponding job in the interval [0,O + (C-1)k + 1]. Hence, we need $log_2(O+(C-1)k+1)$ qubits to get all possibilities of a given job of index $k$ and a sequence of $I\sum_{k} \log_2(O+(C-1)k+1)$ represents a full instance solution. We associate qubit combinations to date by ascending order for each time window, i.e  the combination $\ket{0} \otimes \log_2(O+(C-1)k+1)$  is associated with the earliest time $\min_{i} E_{ik}$ and the highest combination is associated with the date $\max_i(D_{ik})$. Figure \ref{fig:offsetdates} and \ref{fig:offsetQubits} are graphical representations of the allowed dates for the two first jobs for an machine with $O=2$ and $C=4$. Feasible planning corresponds to the path in red, hence labeling feasible solutions involves the combination of qubits associated with the two nodes of a red path.

\begin{figure}[h!]
    \centering
    \begin{minipage}{0.45\textwidth}
        \centering
         \resizebox{\linewidth}{!}{
       \begin{tikzpicture}[shorten >=1pt,->]
  \tikzstyle{vertex}=[circle,fill=black!25,minimum size=17pt,inner sep=0pt]
  
  \node[vertex](G-s) at (0,0) {s};
  
  \foreach \y in {0, 1, 2,3,4,5}
    		\node[vertex] (G-1-\y) at (5,2*\the\numexpr\y-1) {$\y$};
  
  \foreach \y in {0,1,2,3,4,5,6,7,8}
  	\node[vertex] (G-2-\y) at (10,2*\the\numexpr\y-2) {$\y$};

   \draw[red, dashed, thick] ($(G-1-5.north west)+(-0.2,0.2)$) rectangle ($(G-1-2.south east)+(0.2,-0.2)$);
   \node[red, anchor=east] at ($(G-1-5.north west)+(-0.2,0)$) {Dates job 1};
  
   \draw[blue, dashed, thick] ($(G-2-8.north west)+(-0.2,0.2)$) rectangle ($(G-2-2.south east)+(0.2,-0.2)$);
   \node[blue, anchor=east] at ($(G-2-8.north east)+(2.8,0)$) {Dates job 2};

   	\foreach \node in {0,1}
   	{
     		\draw (G-s) -- (G-1-\node);
     		\foreach \othernode in {\the\numexpr\node+3,\the\numexpr\node+2,\the\numexpr\node+1, \the\numexpr\node}
     		{
     			\draw (G-1-\node) -- (G-2-\othernode);
     		}	
    }
    
    \foreach \node in {2,3,4,5}
   	{
     		\draw[red] (G-s) -- (G-1-\node);
     		\foreach \othernode in {\the\numexpr\node+3,\the\numexpr\node+2,\the\numexpr\node+1, \the\numexpr\node}
     		{
     			\draw[red] (G-1-\node) -- (G-2-\othernode);
     		}	
    }

\end{tikzpicture}}
        \caption{Possible job dates \\ machine with an offset of two}
        \label{fig:offsetdates}
    \end{minipage}
    \begin{minipage}{0.45\textwidth}
        \centering
\resizebox{\linewidth}{!}{
\begin{tikzpicture}[shorten >=1pt,->]
  \tikzstyle{vertex}=[circle,fill=black!25,minimum size=17pt,inner sep=0pt]
  
  \node[vertex](G-s) at (0,0) {s};
  
  \foreach \y/\x in {-1/$\ket{000}$, 0/$\ket{001}$, 1/$\ket{010}$, 2/$\ket{011}$,3/$\ket{100}$,4/$\ket{101}$}
  	\node[vertex] (G-1-\y) at (5,2*\the\numexpr\y) {\x};
  \draw[red, dashed, thick] ($(G-1-4.north west)+(-0.2,0.2)$) rectangle ($(G-1-1.south east)+(0.2,-0.2)$);
  \node[red, anchor=east] at ($(G-1-4.north west)+(-0.2,0)$) {Qubit combination};
  \node[red, anchor=east] at ($(G-1-4.north west)+(-0.2,-0.5)$) {job 1};
  \draw[blue, dashed, thick] ($(G-2-8.north west)+(-0.2,0.2)$) rectangle ($(G-2-2.south east)+(0.2,-0.2)$);
  \node[blue, anchor=east] at ($(G-2-8.north east)+(3.5,0)$) {Qubit combination};
  \node[blue, anchor=east] at ($(G-2-8.north east)+(3.5,-0.5)$) {job 2};
  \foreach \y/\x in {-2/$\ket{0000}$,-1/$\ket{0001}$,0/$\ket{0010}$,1/$\ket{0011}$,2/$\ket{0100}$,3/$\ket{0101}$,4/$\ket{0110}$,5/$\ket{0111}$,6/$\ket{1000}$}
  	\node[vertex] (G-2-\y) at (10,2*\the\numexpr\y) {$\x$};

  	\foreach \node in {-1,0}
  	{
    		\draw (G-s) -- (G-1-\node);
    		\foreach \othernode in {\the\numexpr\node+2,\the\numexpr\node+1,\the\numexpr\node, \the\numexpr\node-1}
    		{
    			\draw (G-1-\node) -- (G-2-\othernode);
    		}	
    	}
       	\foreach \node in {1,2,3,4}
  	{
    		\draw[red] (G-s) -- (G-1-\node);
    		\foreach \othernode in {\the\numexpr\node+2,\the\numexpr\node+1,\the\numexpr\node, \the\numexpr\node-1}
    		{
    			\draw[red] (G-1-\node) -- (G-2-\othernode);
    		}	
    	}

\end{tikzpicture}}
        \caption{Qubit combination \\ machine with an offset of two}
        \label{fig:offsetQubits}
    \end{minipage}

\end{figure}

Once we have defined the encoding, the idea is to apply a quantum search algorithm to generate feasible solutions to our planning problem. Textbook quantum search algorithm start by putting the system in a superposition of all possible qubit combinations by applying a Hadamard gate to each qubit. Then, we need to define an oracle circuit which will mark qubit combinations that respect the following constraints:
\begin{itemize}
\item All jobs are within there time windows.
\item Two consecutive jobs of the same machine respect the spacing constraints.
\item Jobs from different machines do not overlap.
\end{itemize}

\subsubsection{Time windows and spacing constraint oracle} \label{subsubsec:FeasiblePathOracle}

As time windows and spacing constraints are local on each machine we can define an oracle which checks whether a given machine planning respect those constraints and combine the result on an unique flag qubit at the end of the circuit. Possible job dates for the first job depend of the machine's offset, hence this circuit takes $\sum_k\log_2(O + Ck + 1)$ qubits to describe job dates and $\log_2(O)$ qubits for the machine offset. Let  $\ket{x^0_i}$ be the binary value associate with the offset of the machine and $(\ket{x^k_i})_{k\geq 1}$ the states associated with all dates of jobs on machine $i$. Remark that the state associated with the $E_i$ corresponds to $\ket{x^1_i} = \ket{O_i}$ and the time windows structure implies that the state associated with $D_i$ correspond to $\ket{x^1_i} = \ket{O_i} + C$. Hence, feasible states for the first job date of the machine correspond to $ 0 \leq x^1_i - O_i \leq C$. Additionally, assumption \ref{ass:Spacing} combined with the association between qubit combination and job dates implies that for a given index $k \geq 1$ $\ket{x^k_i}\ket{x^{k+1}_i}$ corresponds to a feasible solution only if the following inequalities hold : $0 \leq x^{k+1}_i - x^k_i \leq C$. Hence, to check whether the state $\prod_{k \geq 0} x_i^k$ corresponds to a feasible planning with respect to time windows and spacing constraints we need to control that : 
$$ \forall k \geq 0, \qquad 0 \leq x^{k+1}_i - x^{k}_{i} \leq C \label{eq:ineqfeasiblepath}$$

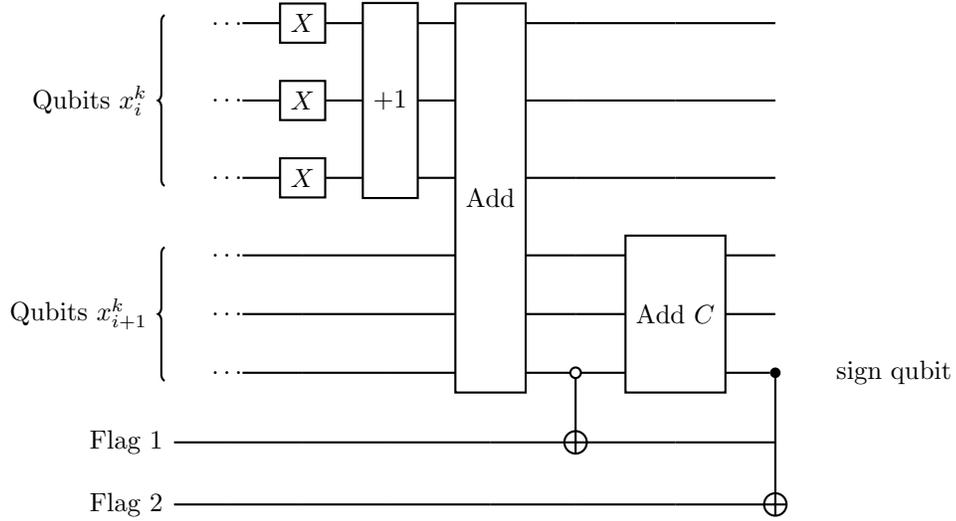
\begin{figure}
    \centering
\begin{quantikz}
    \lstick[3]{Qubits $x^k_i$}     & \ldots & \gate{X} & \gate[3]{+1} & \gate[6]{\text{Add}} &\qw    & \qw                             & \qw  &\\
    				         & \ldots & \gate{X} &                    &                      &   \qw           & \qw                             & \qw  &\\
    				         & \ldots & \gate{X} &                    &                      &     \qw         & \qw                             & \qw  &\\
 \lstick[3]{Qubits $x^k_{i+1}$} & \ldots &  \qw       &              \qw &                      &  \qw             & \gate[3]{\text{Add }C} & \qw   & \\
    					& \ldots &  \qw       &              \qw &                      &  \qw            &                                   & \qw     & \\
    					& \ldots &  \qw       &              \qw &                      & \octrl{1}      &                                   & \ctrl{2} &\rstick{sign qubit}\\
  \lstick{Flag 1} 		  	& \qw   &        \qw  &              \qw &               \qw & \targ{}        & \qw                             &    \qw   &\\
 \lstick{Flag 2} 			& \qw   &        \qw  &              \qw &               \qw & \qw            &\qw                              & \targ{}  & \\
\end{quantikz}
    \caption{The circuit showing the setup of the time windows constraint oracle.}
    \label{fig:circuit_time_windows}
\end{figure}

In arithmetic, negative binary numbers are represented by their ``two's complement'' representation using the help of an extra (qu)bit called ``sign (qu)bit'' such that for any integer number $a$ the sum of $a$ and its two's complement is equal to 0 (because of overflow). Hence, to check if equation \ref{eq:ineqfeasiblepath} holds between two consecutive jobs states $\ket{x^k_i}$ and $\ket{x^{k+1}_{1}}$ we can follow the following steps : 
\begin{enumerate}
\item Take the two's complement of $x^k_i$.
\begin{enumerate}
	\item Invert all qubits by applying the X gate to all of them.
	\item Adding 1 to the number, using the +1 gate (see Box The +1 Gate) in binary.
\end{enumerate}
\item Add $x^k_{i+1}$ and the two's complement of $x^k_i$.
\item Check whether the sign qubit of the resulting $x^k_{i+1}$ is set to zero (e.g., by applying a CNOT onto an extra ancilla ``flag'' qubit).
\item Add the two's complement of $C$.
\item Check whether the sign qubit of the resulting $x^k_{i+1}$ is set to one (e.g., by applying a CNOT onto another ancilla ``flag'' qubit).
\item Undo steps 1, 2 and 4, to restore qubits $x^k_i$ and $x^k_{i+1}$.
\begin{enumerate}
	\item Subtract the two's complement of $C$ from $x^k_{i+1}$.
	\item Subtract $x^k_i$ from $x^k_{i+1}$.
	\item Subtract 1 from $x^k_i$, using the inverse of the +1 gate in binary.
	\item Invert all qubits by applying the X gate to all of them.
\end{enumerate}
\end{enumerate}
The circuit showing steps 1-5 is depicted in figure \ref{fig:circuit_time_windows}.

    \begin{tcolorbox}[title=The +1 Gate,label=+1gate]
        The \textbf{(+1)-gate}, given by $\hat{\mathcal{U}}_{(+1)}$, which constitutes a central building block of the proposed algorithm. The way this gate acts is by taking a state of overall $n$ qubits as input and giving as output the state incremented by one in binary (with overflow): 
        
        \begin{equation}
            \vert \tilde{\textbf{y}}\rangle   \xrightarrow{\hat{\mathcal{U}}_{(+1)}} \vert \tilde{\textbf{y}}\textbf{+1}\rangle  
        \end{equation}
        
        ,where $\tilde{\textbf{y}}=(y_{0}, y_{1}, ..., y_{n-1})$.
        \hfill\break
        \hfill\break
        Some simple examples can be seen below \footnote{recall at this point that that $\textbf{y}_{d}=\sum\limits_{j=0}^{n-1}y_{j}2^{n-(j+1)}$}:
        
        \begin{center}
            \includegraphics[width=0.5\linewidth]{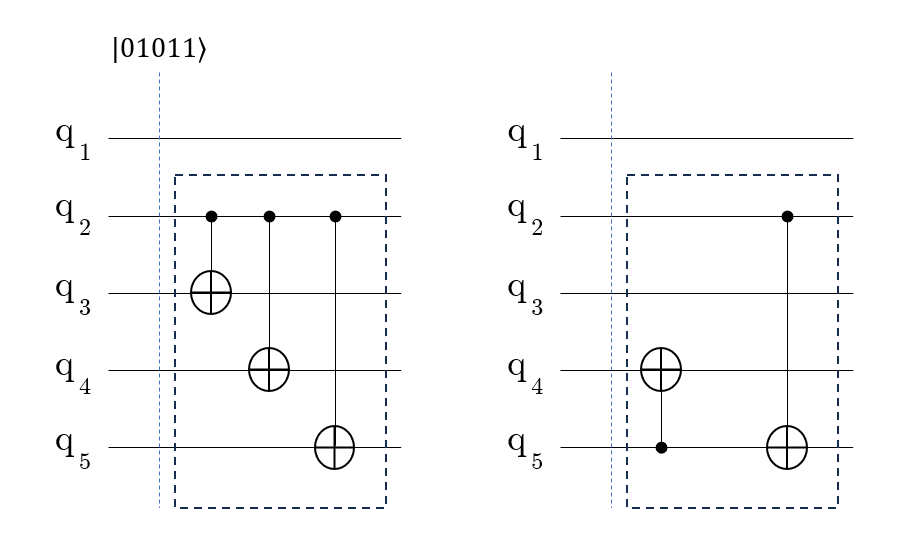}
            \captionof{figure}{Action of $\hat{\mathcal{U}}_{(+1)}$ on a given quantum register}
            \label{fig:enter-label}
        \end{center}

        \begin{align*}
        \ket{01011} & \xrightarrow{\hat{\mathcal{U}}_{(+1)}} \ket{01100} \\
        \frac{1}{\sqrt{2}}\left( \ket{01001} + \ket{11000} \right) & \xrightarrow{\hat{\mathcal{U}}_{(+1)}} \frac{1}{\sqrt{2}}\left( \ket{01010} + \ket{11001} \right) \\
        \end{align*}
        
        There are multiple possible implementations of the above operation. The one used in our case is based on the work \cite{shakeel2020efficient} (defined as the "X-gate" in the article), implemented by a quantum Fourier transform (QFT), followed by a series of single-qubit rotation gates and then followed by another QFT. Implementing the QFTs before and after the rotations requires $\mathcal{O}(n^2)$ gates, for each operation. Here, $n$ is the number of qubits acted upon; as expected, this number which will be changing throughout the algorithm.
        
        \end{tcolorbox}

At this point we have two flag qubits whose states encode the truth value of the statements $0 \leq x^{k+1}_i - x^k_i$ and $x^{k+1}_i - x^k_i \leq C$ respectively. We can assume that they're in the state $\ket{1}$ if the corresponding statement is true. Such pair of flag qubits is defined for any $k$. All of these qubits need to be in the state \ket{1} in order for all of the time-windows and spacing conditions to be fulfilled. Therefore, the oracle should apply a phase $\phi$ to the state of all $\ket{1}$'s and act as identity on all other states. This can be done with a multi-qubit gate implementing the following unitary:
\[
\begin{pmatrix}
1      & 0      & \cdots & 0      & 0 \\
0      & 1      & \cdots & 0      & 0 \\
\vdots & \vdots & \ddots & \vdots & \vdots \\
0      & 0      & \cdots & 1      & 0 \\
0      & 0      & \cdots & 0      & e^{i \phi} \\
\end{pmatrix}.
\]

An alternative (more illustrative) approach would use all of the flag qubits as a control in a multi-controlled NOT gate which acts on yet another ``ultimate'' flag qubit. The state of this flag qubit would then encode the truth value of all the time windows and spacing constraints. The oracle would then just need to apply a phase gate on this one qubit.

This circuit uses $2I(K-1)$ flag ancilla qubits. The depth is dominated by the addition gates, each of which includes a quantum Fourier transform circuit and one layer of 1-qubit gates. That gives depth $\mathcal{O}(n^2)$ where $n$ is the number of qubits the addition gates act on. Since these act on the sets of qubits corresponding to one job, we have $n \approx \log((C-1)K+1)$. The flag qubits can't be set up entirely in parallel (since the qubits $x^k_{i+1}$ used here will be needed in next step to check against qubits $x^k_{i+2}$ and so on ... ), so the depth should be multiplied by a further $(K-1)$ for the number of job pairs to check. However, all machines can be checked in parallel. The final multi-controlled phase gate can be decomposed in a standard way into $\mathcal{O}(I^2)$ two-qubit gates.

\subsubsection{Resource constraint oracle} \label{subsubsec:ResourceConstraintOracle}

Resource constraints bring additionnal limitations for jobs of the same index $k$ from different machines. To define an oracle for this kind of constraint we first need to check separatly for each constraint and each job index if the constraint holds. This requires several ancillary qubits. To avoid confusion, we will divide the ancillary qubits into groups labelled A and B. The information on the qubit number per ancillary qubit group, as well as the number of data qubits in the register, is organized in table \ref{table:ancillasCat} below : 

\begin{table}[h!]
\begin{center}
\begin{tabular}{|l|c|c|}
\hline
   Qubit group & Number of qubits\\
\hline
    $\mathcal{G}_{A}$ & $\frac{I(I-1)}{2}K$ \\ 
    $\mathcal{G}_{B}$ &$1$\\ 
\hline
\end{tabular}
\end{center}
\caption{Group of ancillary qubits and their dimensions}
\label{table:ancillasCat}
\end{table}

For each resource constraint concerning machines $i$ and $i'$, we first check for each job index $k$ if the schedules of the machine do not overlap. This is done by applying CNOT gates controlled by the data qubits $x_{ik}$ acting on the data qubits of $x_{i'k}$. If there is an overlap, this step will set all the data qubits of $x_{i'k}$ to the state $|0\rangle$. 
To verify this, we introduce an ancilla qubit (from the group A) and apply a multi-controlled NOT gate to it controlled by the data qubits $x_{i'k}$ all being in the state $|0\rangle$. Therefore, if the jobs of the machines $i$ and $i'$ overlap, this ancilla qubit is now in the state $\ket{1}$ (and if they don't, it's in the state $\ket{0}$). After this step, we again apply a set of CNOTs controlled by the data qubits $x_{ik}$ acting on the data qubits of $x_{i'k}$, in order to un-compute the first wave of CNOTs (to restore $x_{i'k}$ back to their original state). The entire circuit can be seen in figure \ref{fig:circuit_ancilla_A}.
%
As a last step, we apply a multi-controlled NOT gate, controlled by all the ancillas in group A being in the $|0\rangle$ state and acting on the one ancilla in group B. If this qubit is flipped, it means that no overlaps occure between any two machines in any job index and thus no incompatibility: the state as a whole satisfies the resource constraint. If the qubit is not flipped, the resource constraints are violated for some pair of machines and some job index and this state is now marked for deletion. This part of the circuit is depicted in figure \ref{fig:circuit_ancilla_B}.

\begin{figure}[h]
    \centering
    \begin{minipage}{0.45\linewidth}
        \centering
        \resizebox{\linewidth}{!}{
        \begin{quantikz}
            \lstick[3]{1st p.p. data qubits $x_{ik}$}  & \ldots & \ctrl{3} & \qw & \qw & \qw & \qw & \qw & \ctrl{3}  \\
                                & \ldots & \qw & \ctrl{3} & \qw & \qw & \qw & \ctrl{3} & \qw  \\
                                & \ldots & \qw & \qw & \ctrl{3} & \qw & \ctrl{3} & \qw & \qw  \\
            \lstick[3]{2nd p.p data qubits $x_{i'k}$} & \ldots & \targ{} & \qw & \qw & \octrl{1} & \qw & \qw & \targ{}  \\
                                & \ldots & \qw & \targ{} & \qw & \octrl{1} & \qw & \targ{} & \qw  \\
                                & \ldots & \qw & \qw & \targ{} & \octrl{1} & \targ{} & \qw & \qw  \\
            \lstick{Group A ancilla} 	& \qw & \qw & \qw & \qw & \targ{} & \qw & \qw & \qw 
        \end{quantikz}}
        \caption{The circuit showing how each ancilla qubit in group A is set up using the data qubits $x_{ik}$ and $x_{i'k}$.}
        \label{fig:circuit_ancilla_A}
    \end{minipage}
    \hfill
    \begin{minipage}{0.45\linewidth}
        \centering
        \resizebox{\linewidth}{!}{
        \begin{quantikz}
            \lstick[3]{Group A ancillas}  & \ldots & \octrl{1} & \qw  \\
                                & \ldots  & \octrl{1} & \qw \\
                                & \ldots  & \octrl{1}  & \qw \\
            \lstick{Group B ancilla} 	& \qw &  \targ{} & \qw
        \end{quantikz}}
        \caption{The quantum circuit showing how ancilla B is acted on controlled by ancillas A.}
        \label{fig:circuit_ancilla_B}
    \end{minipage}
\end{figure}

From the above we can see that the qubit requirements only grow quadratically with the number of machines and slightly more than linearly with the number of time steps, namely $\mathcal{O}(I^{2}K)$.

The comparators themselves have depth ~$\mathcal{O}(1)$. If we use many A-ancillas, we can apply all of these in parallel; however if we reuse A-ancillas, we need to apply these in sequence $\frac{I(I-1)}{2}K$ times. The multi-controlled X controlled by the A qubits acting on the B-qubits itself will need circuit of depth $\mathcal{O}(\log((C-1)K+1)$ \cite{saeedi2013linear} and we need to apply it $\frac{I(I-1)}{2}K$ times. Again, depending on whether we recycle A-ancillas, we might need to apply them in series as opposed to in parallel. After that we only need to apply one (inversed) multi-controlled X, which takes $\mathcal{O}(\frac{I(I-1)}{2}K)$ elementary gates. With this, the leading term in the number of gates is either $\mathcal{O}(\frac{I(I-1)}{2}K)$ if the A ancillas were paralellized or $\mathcal{O}\left[\frac{I(I-1)}{2}K\log(CK)\right]$ if they were reused. 



\subsubsection{Amplitude amplification}

The amplitude amplification algorithm selectively amplifies the the states with the flag qubit set to $\ket{1}$ and decrease the other states amplitude. The standard amplitude amplification algorithm consists of $N/M$ steps of applying the following two machineary operators:
\begin{align}
    S_P &= 1 - 2 \ket{0}\bra{0}_{\tilde{Q}} \\
    S_\phi &= 1- 2\ket{\phi}\bra{\phi} 
\end{align}

Here $\ket{\phi}$ is the state of the system before starting the amplitude amplification algorithm. The unitary $S_P$ is easy to construct. It's just the Z-gate acting on the flag ancilla as a flip. The unitary $S_\phi$ is more complicated. In particular, its function is to flip the phase of the state $\ket{\phi}$ while keeping the phase of all other states the same. This is implemented by initially "un-computing" the entire circuit (so that the state $\ket{\phi}$ is transformed back to the source state $\ket{000(...)}$ again. This is performed by the overall conjugate operation (application of each of the operations applied starting from the source state, conjugated and in inverse order) on the full quantum register ($\mathcal{Q}$), followed by a multi-anti-controlled Z-gate: in such way, the phase changes if and only only if all qubits are in the state $\ket{0}$. Finally, we apply the whole circuit again. This obviously requires to calculate the entire circuit \textit{twice}. The complete operation is found below:

\begin{equation}
    S_{\phi}\vert \phi\rangle=\Big(\mathcal{U} \cdot Z \cdot\mathcal{U}^{\dagger}\Big)\vert\phi\rangle
\end{equation}

, where $\mathcal{U}$ and $\mathcal{U}^{\dagger}$, correspond to the full set of operations throughout the whole circuit and their inverse ones, respectively and with $\mathcal{U}^{\dagger}\vert\phi\rangle\to\vert 000..\rangle$.


Standard amplitude amplification is usable only if we know \textit{accurately} the proportion of qubit combinations that satisfy the constraints implemented in our oracles. Here, as this information in unknown, we employ "fixed-point amplitude amplification," which can be employed regardless of the number of the correct combinations \cite{gilyen2019quantum}.In particular, if there are $M$ correct states among a superposition of $N$ basis states, standard amplitude amplification needs to be used with \emph{exactly} $\sqrt{N/M}$ iterations (queries to the oracles) to find them. Fixed-point amplitude amplification needs \emph{at least} $\sqrt{N/M}$ iterations to find the correct states. Even if $M/N$ is unknown, one can use the fixed-point amplitude amplification combined with  \emph{exponential search}\cite{Bentley1976almost} to find the solutions with only small constant overhead. The seach therefore adds a factor of $\mathcal{O}\left(\frac{N}{M}\right)$ to the circuit depth (multiplying the depth of the two oracles).

\subsection{Space size analysis} \label{subsec:FSAnalyse}

\begin{figure}[h]
    \centering
    \begin{minipage}{0.45\textwidth}
        \centering
        \includegraphics[width=\textwidth]{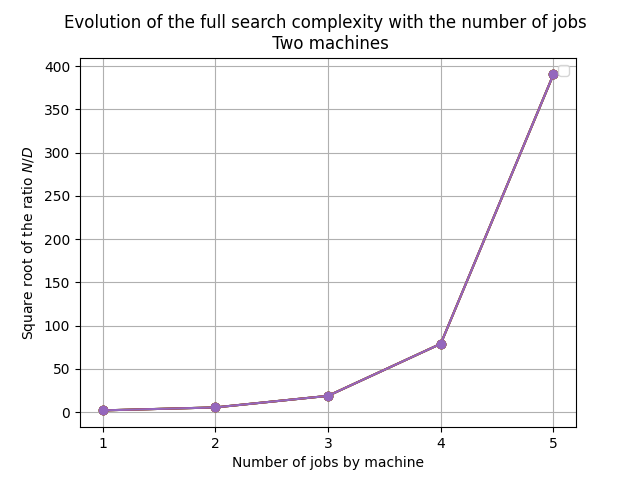}
        \caption{Full Search Complexity for Two machines}
        \label{fig:FSC_2units}
    \end{minipage}
    \hfill
    \begin{minipage}{0.45\textwidth}
        \centering
        \includegraphics[width=\textwidth]{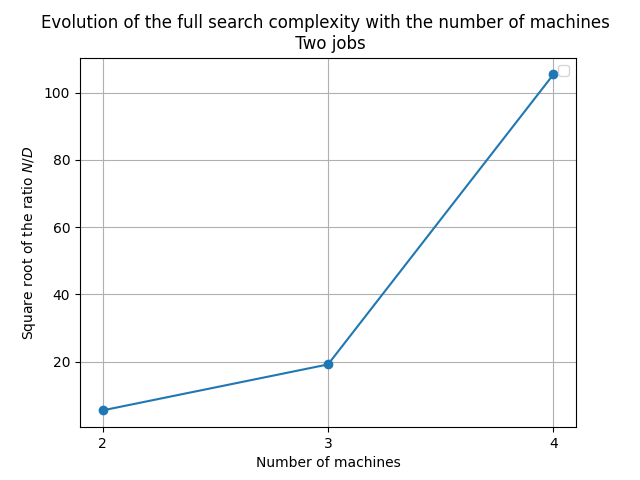}
        \caption{Full Search Complexity Number of machines Impact}
        \label{fig:FSC_WithUnits}
    \end{minipage}
\end{figure}

As mentioned in Section \ref{sec:QSearch}, quantum search demonstrates a theoretical quadratic speedup for searching in an unsorted space. However, the complexity of the algorithm depends on the ratio between the total number of possible qubit combinations and the number of marked elements. For our simplified use case, in the case without offset, the full search space is of size $2^{I.\sum_{k} \log_2(C + (C-1)k)}$. The marked elements correspond to combinations of paths which do not violate the resource constraint. Figures \ref{fig:FSC_2units} gives the evolution of the square root of the ratio between the total number of qubit combinations and the number of solutions (obtained by brute force approach in this example) with the number of jobs for two machines. Figure \ref{fig:FSC_WithUnits} displays the same evolution for different number of machines $I$. 

We observe that the expected number of iterations increases exponentialy both with the number of jobs (figure \ref{fig:FSC_2units}) and with the number of machines (figure \ref{fig:FSC_WithUnits}), making the search impossible in practice for large instances. In the next section we propose to an improve quantum search approach by building an initial state superposition of smaller size, namely corresponding to the qubit combination which respect the time windows and spacing constraints, to decrease the ratio $D/N$ and allow potential futur applications of quantum search.

\section{Quantum walk inspired state-space reduction exploiting regular state-graph structure}  \label{sec:SSR} 

As detailed in the previous section quantum search algorithms start from an initial state (superposition of states) supplemented by an oracle which marks ``good'' elements and iteratively rotate from this initial state to superpositions where marked element has a high density. As the number of steps required for quantum search algorithms depends on the ratio between good and bad elements, starting from an initial state where this ratio is as high as possible allow to perform less iteration to ensure measuring a good elements.
In this section, we introduce a way to increase the ratio between ``good'' and ``bad'' elements in the initial state superposition when the solution space induced by a subset of constraints corresponds to paths in an oriented rooted graph. This method is inspired from quantum walk approaches described in the following section. 

\subsection{Introduction to quantum walk}

A quantum walk is a straightforward generalization of a classical random walk, consisting of discrete steps. In each step, the “walker” randomly chooses a step to take from a set of possible steps in a structured space, which in our case corresponds to the rooted directed graph shown in Figure \ref{fig:GraphExemple}. At each step, the walker is located at one of the nodes of the graph, and the set of possible steps corresponds to the outgoing edges of that node. An edge is randomly selected, and the walker moves along the edge to the neighboring node. Since the position of the walker is random, it can be expressed as a probability distribution.

In a quantum random walk (QRW) \cite{kempe2003quantum}, instead of randomly choosing a step, the walker takes a superposition of all available steps. Therefore, the position of the walker is not represented by a probability distribution, but by a superposition. A quantum walk is described by a quantum state that evolves under both a unitary operator (representing the “coin flip” step) and a conditional shift operator (representing the walker’s movement). While classical random walks use probabilities to describe the likelihood of transitions between states, QRWs use quantum amplitudes to govern the spread of the walker’s quantum state over the space. This kind of approach comprises the following component :

\begin{itemize}
    \item Coin Flip Operator ($\mathcal{C}$): This operator acts on a separate "coin" space ($\mathcal{H}_{C}$) associated with each vertex of the graph. It represents the probabilistic decision the walker makes at each step. Formally, it's a unitary operator acting on the coin space. If we have 
    N vertices, the coin space is typically represented by a 
    $N \times N$ unitary matrix.
    
    \item Shift Operator ($\mathcal{S}$): This operator describes the movement of the walker conditioned on the outcome of the coin flip. It determines how the walker transitions between vertices based on the state of the coin. Formally, it's a conditional shift operator that depends on the coin state.
    \end{itemize}
    
    Consequently the hilbert space of the system is given by the product of the coin and shift hilbert spaces: 
    \begin{equation}
        \mathcal{H}= \mathcal{H_{C}}\otimes \mathcal{H_{S}}
    \end{equation}
    The evolution of the walker's state from time $t$ to time $t+1$ is given by the following operation:
    
    \begin{equation}
        \vert\psi(t+1)\rangle= \mathcal{S}\Big(\mathcal{C}(\vert\psi(t)\otimes\vert 0\rangle)\Big)
    \end{equation}
    
    This process can be iterated for multiple steps to simulate the evolution of the walker over time. In the present study, this aspect is emphasized because quantum walk search algorithms are appealing for their at most quadratic advantage over classical algorithms in exploring graphs. This advantage ultimately stems from the same principles that provide the speed-up in Grover’s algorithm.

\subsection{Quantum walk inspired state space reduction }

\subsubsection{Idea of the state space reduction}

As discussed in the previous sections, quantum random walks (QRWs) facilitate the construction of a superposition of states within a graph-like structure, specifically a tree in our case. Our work exploit this technic to build a superposition of paths within the network that correspond to feasible schedules for a set of machines, subject to operational constraints such as time windows and spacing. This is accomplished by allowing the QRW to explore the corresponding space of paths on the tree, with each new step performed on new qubit data. 

\vspace{0.2cm}

This algorithm, called \textit{Feasible path algorithm}, will act as follows : 

\begin{enumerate} 
    \item Each set of data qubits represents the state associated with a given machine. The walk starts by the inititialization of the system with $\ket{x_i^0}$ initialized to the offset of the machine and all other data qubits in $\ket{0}$ state:  
    \begin{equation*}
    \ket{s_ i}= \ket{x^0_i}\otimes\prod_{k=1}^{K}\ket{x_i^k}= \ket{x^0_i}\otimes\ket{00..0}
    \end{equation*}

    Now, each step of the quantum walk corresponds to a decision or transition between possible scheduling choices of a single or multiple machines.

    \item We iteratively apply the quantum walk-based scheme on all qubits associated with jobs $k$ $\ket{x^k_i}$ and store the result in qubits $\ket{x^{k+1}_i}$. This consists, to the application of the coin and shift operators: 
    \begin{enumerate}
        \item   The coin operator ($\mathcal{C}$) is applied introducing superposition over different possible paths (schedules). This represents the various scheduling choices and its dimension corresponds to the number of possible paths: dim$[\mathcal{C}]=n_{\text{paths}}$
	\hfill\break
        \item The shift operator ($\mathcal{S}$) develops the walk through the space of paths by extending the partial schedule to the next job index.
    \end{enumerate}
    After $k^{th}$ steps of the algorithm the configuration is the following: 
    \begin{equation*}
        \ket{s_k}= \prod_{k' \leq k}\ket{x_i^k}\otimes \prod_{k'=k+1}^{K}\ket{0} = \prod_{k' \leq k}\Big(\mathcal{S}_{[i]}\mathcal{C}_{[i]}\Big)^{k}\ket{x_i^0}\otimes \prod_{k'=k+1}^{K}\ket{0}
    \end{equation*}
    where $\prod_{k' \leq k}\ket{x_i^k}$ correspond to the superposition of path from the source to all node of depth $k$ in the graph.
    \hfill\break
    \item After $K$ steps the walker as explored all path of the graph and hence the data qubit are in superposition of all possible path :
    \begin{equation*}   
        \ket{s_k} = \prod_{k' \leq k}\Big(\mathcal{S}_{[i]}\mathcal{C}_{[i]}\Big)^{K}\ket{x_i^0}
    \end{equation*}
    At this point, we have obtained a superposition of all possible paths, which can be used as an an initial state for the amplitude amplification algorithm with the resource constraint oracle described in section \ref{subsubsec:ResourceConstraintOracle}.
\end{enumerate}

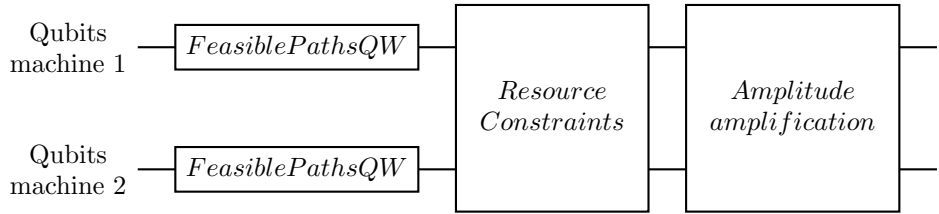
\begin{figure}[ht!]
     \label{fig:alg_circ}
     \begin{quantikz}
     \lstick{Qubits \\ machine 1}  
      & \gate{Feasible Paths QW} & \gate[wires=2]{\begin{array}{c}Resource \\ Constraints \end{array}} & \gate[wires=2]{\begin{array}{c}Amplitude \\ amplification \end{array}} 
     & \rstick{}\qw \\
     \lstick{Qubits \\ machine 2}  &  
      \gate{Feasible Paths QW} & &
     & \rstick{} \qw \\
     \end{quantikz}
     \caption{\small{Structure of the reduced search algorithm}}
     \label{fig:RSAlgo}
 \end{figure}

The method described above constitutes the initial phase of our approach to achieving a potentially enhanced quantum search speedup by integrating two techniques. First, we employ a quantum walk to explore the space of paths, introducing a path superposition that serves as the initial state for the quantum search. Notably, this new algorithm requires only the oracle for resource constraints and not the feasible path oracle anymore as the initial state includes only elements that satisfy the time windows and spacing constraints. Figure \ref{fig:RSAlgo} show the structure of the improve algorithm where \textit{FeasiblePathsQW} represents the initial state construction.

\subsubsection{Practical implementation of the quantum-walk inspired scheme} \label{subsubsec:QW}

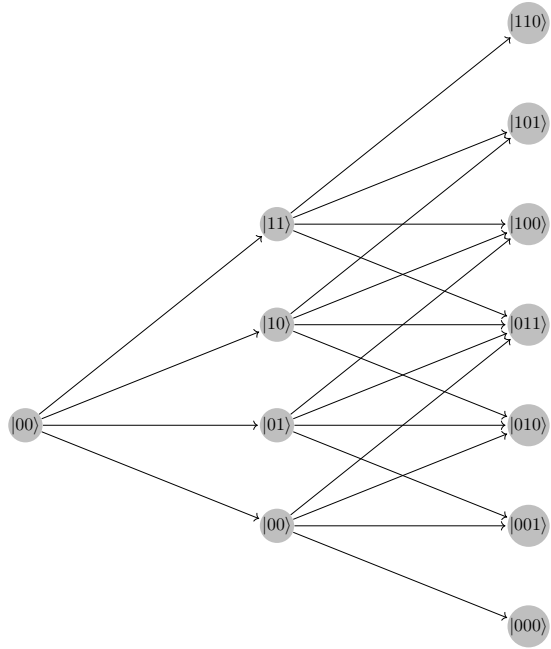
\begin{figure}[!ht]
    \centering
    \resizebox{.6\linewidth}{!}{
    \begin{tikzpicture}[shorten >=1pt,->]
      \tikzstyle{vertex}=[circle,fill=black!25,minimum size=17pt,inner sep=0pt]
      
      \node[vertex](G-s) at (0,0) {$\ket{00}$};
      
      \foreach \y/\x in {-1/$\ket{00}$, 0/$\ket{01}$, 1/$\ket{10}$, 2/$\ket{11}$}
                \node[vertex] (G-1-\y) at (5,2*\the\numexpr\y) {$\x$};
      
      \foreach \y/\x in {-2/$\ket{000}$,-1/$\ket{001}$,0/$\ket{010}$,1/$\ket{011}$,2/$\ket{100}$,3/$\ket{101}$,4/$\ket{110}$}
          \node[vertex] (G-2-\y) at (10,2*\the\numexpr\y) {$\x$};

          \foreach \node in {-1,0,1,2}
          {
                \draw (G-s) -- (G-1-\node);
                \foreach \othernode in {\the\numexpr\node+2,\the\numexpr\node+1,\the\numexpr\node, \the\numexpr\node-1}
                {
                    \draw (G-1-\node) -- (G-2-\othernode);
                }	
            }

    \end{tikzpicture}}
    \caption{Instance for the feasible path demonstration}\label{fig:GraphExemple2}
\end{figure}

This section describes the practical implementation of one iteration of the second step of the algorithm, specifically how we apply a quantum walk-inspired scheme at each iteration. For simplicity, we will illustrate this process using the example of two jobs with an offset of zero as represented in figure \ref{fig:GraphExemple2}.

According to this algorithm, every step is performed onto different sets of qubits. Note that that each node of the first job can always reach the node where the two last qubit are identical to his state and the three node corresponding respectively to add (+1), (+2) and (+3) to the value of this node. Hence, our QW is based on the (+1)-gate described in box \ref{+1gate}.

\vspace{0.5cm}

We now have every building block we need for our feasible path algorithm. This circuit use two ancilla qubits $q_1^\ast$ and $q_2^\ast$ that we will call  \textit{coin ancilla qubits} to underline the similarity with quantum walk. One step of this circuit is discribed as : 

\begin{enumerate}
    \item We copy the state of the data qubits $\ket{x^k_i}$ into the data qubits of the following $\ket{x^{k+1}_i}$ using a serie of CNOT gates. 
    \item We apply the Hadamard gate to the coin ancilla qubits  $q_1^\ast$ and $q_2^\ast$ to put them in a superposition of the four computational basis states.
    \item We apply the (+1)gate on the qubits $\ket{x^{k+1}_i}$ by the state of the $q_1^\ast$ ancilla qubit.
    \item We apply the +1 gate twice on the qubits $\ket{x^{k+1}_i}$ controlled by the state of the $q_2^\ast$ ancilla qubit. 
\end{enumerate}

Let us now describe the application of this circuit to the construction of the states associated with graph \ref{fig:GraphExemple2}. As the first job is represented by the first 2 qubits and all possible combinations are allowed we do not need to use our circuit on them but we only apply a Hadamard operation, putting them in equal superposition:
$$
 \vert q_1 q_2 \rangle \equiv \vert 0\rangle \otimes \vert 0\rangle  \xrightarrow{\text{\textit{H}}}{} \frac{1}{2}\Big(\vert00\rangle+\vert01\rangle+\vert10\rangle+\vert11\rangle\Big).
$$
This represents all the possible paths we can take from the starting (source) node, namely the first step of the quantum walk. Let's now apply our circuit to build the full state superposition of all path leading to the possible job dates for the second job :
\begin{enumerate}
    \item We apply two CNOT gates controlled by $q_1$ and $q_2$ acting on $q_4$ and $q_5$ respectively. Now the state of the first 5 qubits is the following ($q_3$ is yet intact, so it remains in the state $\ket{0}$:
    $$
    \ket{q_1 q_2 q_3 q_4 q_5} = \frac{1}{2} \left( \ket{00000} + \ket{01001} + \ket{10010} + \ket{11011}  \right)
    $$    
    \item We put $q_1^\ast$ and $q_2^\ast$ in superposition using hadamard gates.
    \item We apply the (+1)-gate on the qubits $q_2, q_3$ and $q_4$ controlled by the state of the $q_1^\ast$ ancilla qubit. We currently have the following superposition : 
    \begin{align*}
        \ket{q_1q_2\underline{q_3q_4q_5}}\ket{q_1^\ast} = \frac{1}{2\sqrt{2}} \big( & \ket{00\underline{000}}\ket{0} + \ket{01\underline{001}}\ket{0} + \ket{10\underline{010}} \ket{0} +  \ket{11\underline{011}}\ket{0} \\
        & \ket{00\underline{001}}\ket{1} + \ket{01\underline{010}}\ket{1} + \ket{01\underline{011}} \ket{1} +  \ket{01\underline{100}}\ket{1} \big)
        \end{align*}
    \item We apply the (+1)-gate twice on the qubits $q_2, q_3$ and $q_4$ controlled by the state of the $q_2^\ast$ ancilla qubit. We obtain :  
    \begin{align*}
    \ket{q_1q_2\underline{q_3q_4q_5}}\ket{q_1^\ast q_2^\ast} = \frac{1}{4} \big( & \ket{00\underline{000}}\ket{00} + \ket{00\underline{001}}\ket{01} + \ket{00\underline{010}} \ket{10} +  \ket{00\underline{011}}\ket{11} \\
    & \ket{01\underline{001}}\ket{00} + \ket{01\underline{010}}\ket{01} + \ket{01\underline{011}} \ket{10} +  \ket{01\underline{100}}\ket{11} \\
    & \ket{10\underline{010}}\ket{00} + \ket{10\underline{011}}\ket{01} + \ket{10\underline{100}} \ket{10} +  \ket{10\underline{101}}\ket{11} \\
    & \ket{11\underline{011}}\ket{00} + \ket{11\underline{100}}\ket{01} + \ket{11\underline{101}} \ket{10} +  \ket{11\underline{110}}\ket{11} \big)
    \end{align*}
\end{enumerate}

Here each row corresponds to one path at the first crossroads (the state of the first 2 qubits). Each column corresponds to one path at the second crossroads. The ancillary qubits identify this choice, but the state of the qubits $q_3, q_4$ and $q_5$ also depends on the state of the first two qubits.

\hfill\break

Overall, to construct the feasible paths, we require only $\log_2{C}$ ancilla qubits to use them as coin qubits. 




\begin{equation}
  \mathcal{D}_{(q_{d}, q^{\*})} = \sum\limits_{k=1}^K \mathcal{O}(\log((C-1)K+1)^2) \sim \mathcal{O}(K \log(CK)^2)
\end{equation}


\subsection{Reduced search space size analysis} \label{subsec:RSAnalyse}

\begin{figure}[ht!]
    \centering
    \begin{minipage}{0.45\textwidth}
        \centering
        \includegraphics[width=\textwidth]{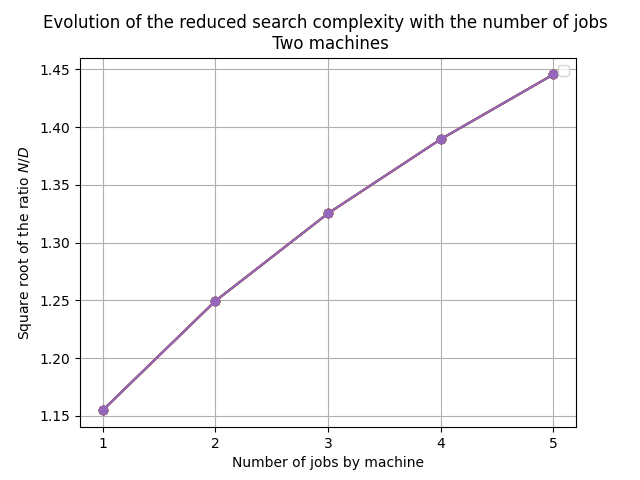}
        \caption{}
        \label{fig:RSC_2units}
    \end{minipage}
    \hfill
    \begin{minipage}{0.45\textwidth}
        \centering
        \includegraphics[width=\textwidth]{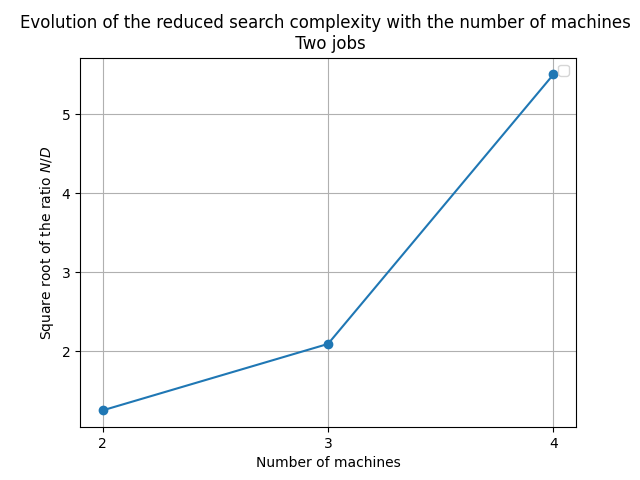}
        \caption{Full Search Complexity Number of machines Impact}
        \label{fig:RSC_WithUnits}
    \end{minipage}
\end{figure}

As in section \ref{subsec:FSAnalyse}, we analyzed the square root ratio between the number of marked elements and the number of possibilities in the reduced search space. Here, the reduced search space is the product of all combinations of possible paths for each machine, resulting in a size of $C^{KI}$. The number of marked elements remains as computed by our brute force approach. Once again, figure \ref{fig:RSC_2units} presents the evolution of the ratio with the number of jobs for two machines, and figure \ref{fig:RSC_WithUnits} the impact of the number of machines. The increase in complexity is quasi-linear with the number of jobs demonstrating the high potential of the approach compared to the exponential ratio of the full search algorithm. The complexity still increase exponentialy with the number of machines but starting from 2 for two machines it ends arround 6 for 4 machines while it was over 100 for the full search space (figure \ref{fig:FSC_WithUnits}). 

\section{Numerical results} \label{sec:Results}

\subsection{Quantum search implementation}

To demonstrate our algorithm, we build the corresponding quantum circuit in Qiskit 0.45.0 \cite{javadi2024quantum}. At the high level, the circuit is relatively simple:
\begin{enumerate}
\item From the problem description, calculate how many qubits are needed and initialize them (starting in \ket{0})
\item Put the qubits in the initial state for the search.
\item Apply the QSVT-based fixed point search.
\item Optionally, measure the state and check that the search was successful.
\end{enumerate}

The main part of the circuit is the QSVT-based fixed point search. In essence, the fixed point search is alternating application of the two oracles with a custom list of phase angles, generated via the \texttt{pyqsp} Python package  \cite{martyn2021grand, chao2020finding, dong2021efficient}, applied to the marked states (not merely flipping the phase). 

There are two sets of oracles: one for the naive quantum search described in section \ref{sec:QSearch} and the other for state-space-reduction approach presented in section \ref{sec:SSR}. In both cases one oracle marks the initial state and the other oracle marks the good states.
In the state-space-reduction approach, the initial state is built by a quantum-walk-like approach. Two qubits are used as a 4-state coin: in each time step, the state of qubits from the previous time step is coherently copied onto the current-time-step register and controlled by the coin qubits, it is incremented by +0/+1/+2/+3 in binary (using quantum arithmetic based on Fourier transform). Thus an initial state is created, which consists of the coherent superposition of all feasible job paths of all power-plants. An oracle to mark this state is built by un-computing the initial state back to all-\ket{0} state, marking it and then re-computing the initial state back. The second oracle is then needed only to check for overlaps of jobs between machines (resource constraints), which is done by comparing the qubit state of the corresponding pairs of registers (using CNOT gates).

For the approach without the state space reduction, the initial state is the standard superposition of all computational basis states (obtained by applying Hadamard gates to all the data qubits). The first oracle marks this state (by "uncomputing" the Hadamard gates and then marking the all \ket{0} state). The second oracle then has to check for both the feasibility constraints and the resource constraints. The resource constraints are checked the same way as before (comparing qubit states by CNOT gates), but on top of that the oracle also checks the feasibility of the job paths for each machine. This is done again by using quantum arithmetic (based on Fourier transform) to compare each pair of consecutive jobs of each machines. If the difference between the two consecutive jobs is outside a pre-defined range (in our case 0-3), the state is marked as bad.

The fixed-point search algorithm needs to have a pre-defined number of iterations (which is also necessary to calculate the angles in \texttt{pyqsp}), but as long as this is large enough (see sec. \ref{sec:SOTA}), we expect the final state to contain predominantly correct solutions. While there is already quantum hardware with sufficient number of qubits to demonstrate a toy-sized version of our problem, the limiting factor is circuit depth. Repeated application of the QSVT oracles requires either error-corrected quantum HW (not yet available) or perfect quantum simulators, which is what we use in our experiment. The test has been performed on AMD Ryzen 5 PRO 5650U with Radeon Graphics 2.30 GHz and 16.0 Go of RAM. 

\subsection{Impact of the search space reduction}

In this section, we analyze the numerical improvement induced by the state space reduction circuit. As presented in \cite{yoder2014fixed}, the fixed-point quantum search does not exhibit a monotonic improvement in the probability of finding a marked element. Instead, it has two phases: the probability first increases until it reaches 100\%, then follows a pseudo-sinusoidal pattern with a minimum that depends on the settings of the fixed-point search.

\hspace{0.5cm}

\begin{figure}[ht!]       
    \centering
    \includegraphics[width=0.8\linewidth]{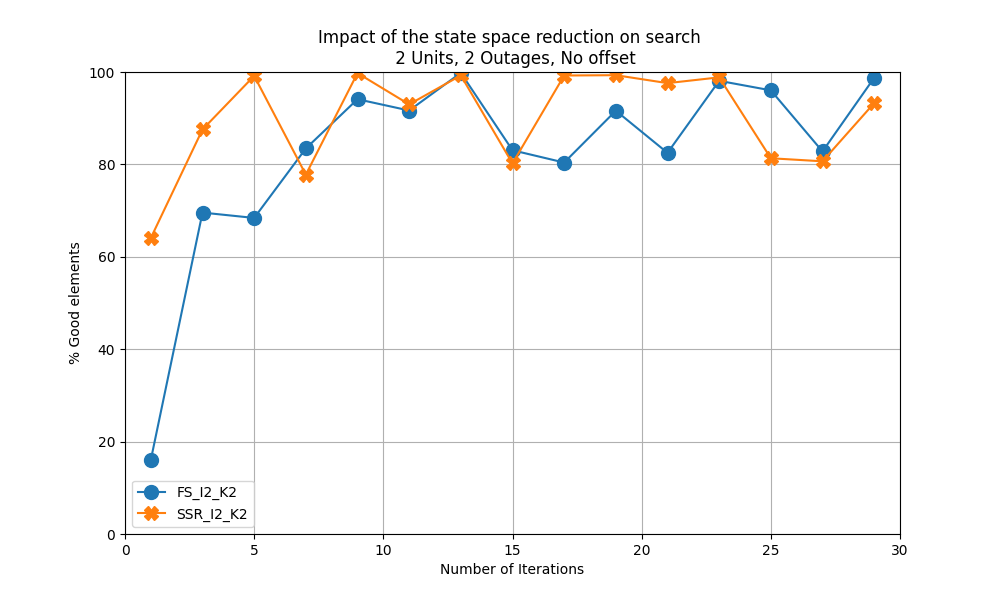}      
    \caption{State space reduction : 2 machines 2 jobs}
    \label{fig:SSR_I2_K2}
\end{figure}

Figure \ref{fig:SSR_I2_K2} compares the percentage of marked elements obtained depending on the number of iterations of the fixed-point quantum search for both the full search ($FS_I2_K2$) and the reduced search ($SSR_I2_K2$), using two machines and two jobs without offsets. As expected, both methods show an initial increasing phase followed by the pseudo-sinusoidal pattern. The percentage of marked elements is higher than 60\% at the beginning for the reduced search and below 20\% for the full search. The reduced search reaches a percentage of 99.16\% of marked elements after only 5 iterations and then enters the pseudo-sinusoidal phase, with the percentage of marked elements never falling below 77\% in subsequent iterations. In contrast, the full search requires 13 iterations to reach a percentage of marked elements greater than 99\% and then enters the pseudo-sinusoidal phase. This result confirm the theoretical improvement induced by the state space reduction.

\begin{figure}[ht!]
    \centering
    \begin{minipage}{0.45\textwidth}
        \centering
        \includegraphics[width=\textwidth]{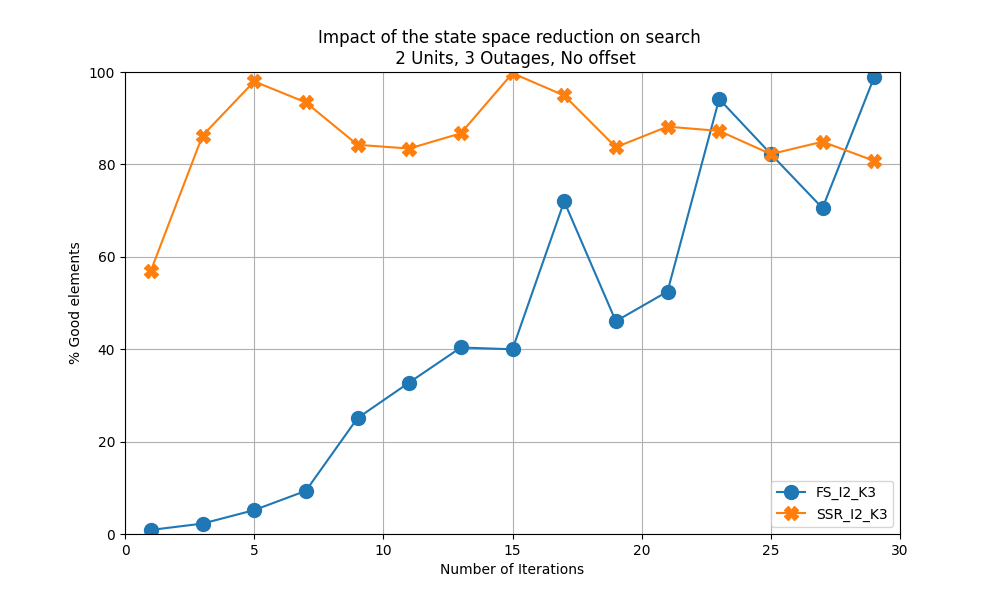}
        \caption{State space reduction : Impact of the number of job $K$}
        \label{fig:SSR_I2_K3}
    \end{minipage}
    \hfill
    \begin{minipage}{0.45\textwidth}
        \centering
        \includegraphics[width=\linewidth]{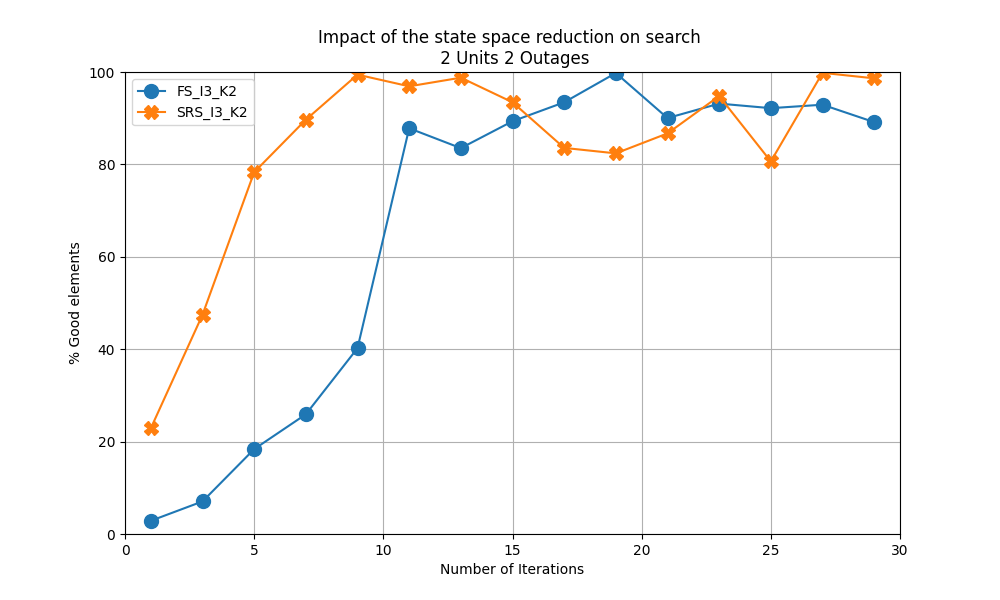}
        \caption{State space reduction : Impact of the number of machines $I$}  
        \label{fig:SSR_I3_K2}    
    \end{minipage}
\end{figure}

Figures \ref{fig:SSR_I2_K3} and \ref{fig:SSR_I3_K2} respectively show the results for 2 machines with 3 jobs and 3 machines with 2 jobs using our approaches. Due to emulator limitations, we could not extend beyond these values. As discussed in sections \ref{subsec:FSAnalyse} and \ref{subsec:RSAnalyse}, the square root ratio between the number of possibilities and the number of solutions increases linearly with the number of jobs in the reduced search, while it increases exponentially in the full search. This is clearly illustrated in figure \ref{fig:SSR_I2_K3}, where the full search space starts with a very low probability of finding good elements and requires 29 iterations to reach a probability higher than 99\%. In contrast, the reduced search starts at approximately 57\% and still reaches 99\% in just 5 iterations.

The impact of increasing the number of machines (figure \ref{fig:SSR_I3_K2}) is more challenging for the reduced search approach, as the search complexity increases exponentially. Here, a high probability of finding a solution is achieved after 9 iterations, which is still better than the full search approach, which takes twice this number of iterations to reach this value.

\newpage

\subsection{Potential scaling up of the approach}

Our implementation serves as a proof of concept for state-space reduction techniques applied to a simplified scheduling problem. However, as mentioned in Section \ref{sec:IndustrialUseCase}, it corresponds to a simplified model for outage planning problems of production units. To apply our approach to real-size industrial instances, we need to scale it up, which will require a number of qubits and a circuit depth that we currently cannot achieve on a quantum emulator. Here, we analyze these future requirements for our approach. Note, however, that our circuit has not been designed to optimize circuit depth, which means it might be possible to reduce this depth by dedicated work on this aspect. We also assume that we can reset ancilla qubits during computation, which decrease the number of required qubits but increase the circuit depth. \\

An industrial instance will typically have the following range of values for the problem parameters : 
\begin{itemize}
    \item Number of machines (units): $I \in [2,100]$
    \item Number of jobs (outage): $K \in [3,10]$
    \item Number of possibilities per time steps: $C \in [10,20]$  
    \item The offset $O \in [5,20]$
\end{itemize}

In the following of this section we will fix $C=15$ and $O=10$ as these two factors have only a slight impact on the qubit and depth requirements.

\subsubsection{Number of qubits}

The amount of qubit required by the full search space approach is : 
\begin{itemize}
\item Data qubits :   $I\sum_{k}\log(O + (C-1)K+1)$
\item Ancillas qubits for resource constraints : $\frac{I(I-1)K}{2}+1$ 
\item Ancillas qubits for time constraints oracle : $2I(K-1)$
\end{itemize}

Overall, for high values of I and K the number of ancilla qubits is dominated by the resource constraint one which gives the following formula for the amount of qubits:
\begin{equation}
I\sum_{k}\log((C-1)K+1 ) + \frac{I(I-1)K}{2}+1 \label{eq:qubitrequire}
\end{equation}

\begin{figure}[ht!]
    \centering
    \begin{minipage}{0.45\textwidth}
        \centering
        \includegraphics[width=\textwidth]{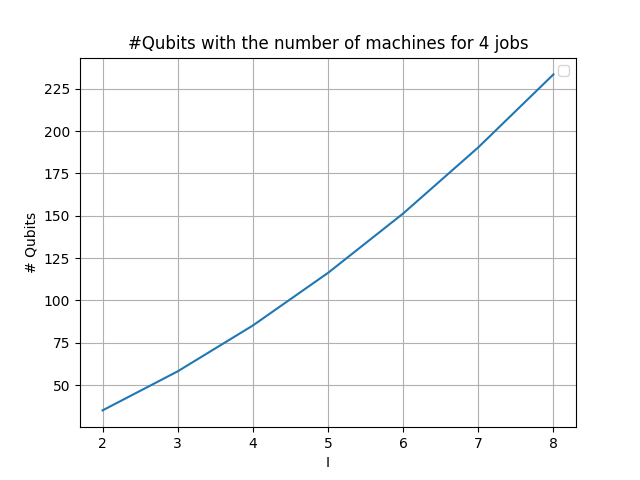}
        \caption{Qubit requirement industrial instances 2 to 6 machines}
        \label{fig:Qubits2-6}
    \end{minipage}
    \hfill
    \begin{minipage}{0.45\textwidth}
        \centering
        \includegraphics[width=\linewidth]{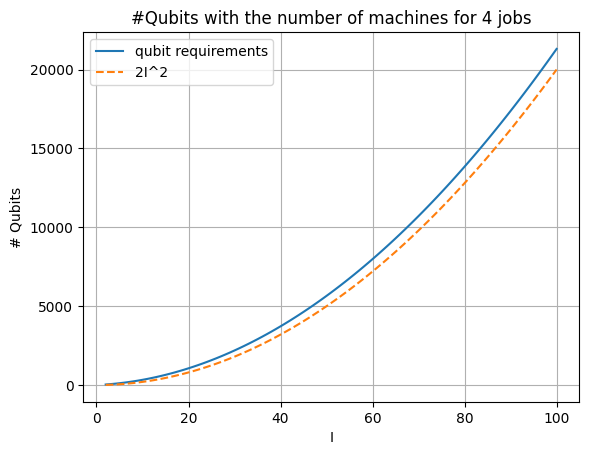}
        \caption{Qubit requirement industrial instances 2 to 100 machines}  
        \label{fig:Qubits2-100}    
    \end{minipage}
\end{figure}

The state-space reduction-based approach requires only $\log_2(C)$ ancilla qubits to perform the state-space reduction. These qubits can also be utilized for the resource constraint oracle. Consequently, the state-space reduction does not necessitate any additional qubits. Figure \ref{fig:Qubits2-6} illustrates the qubit requirements for the aforementioned setting with four jobs and a number of machines ranging from 2 to 6.

It is observed that 100 qubits are sufficient to solve instances with up to 4 machines, which is the typical size for a single production geographical site. Two sites might be accounted for with 230 qubits. Figure \ref{fig:Qubits2-100} extends the analysis to 100 machines and compares it with the curve associated with $2I^2$, demonstrating that this curve serves as a good approximation for the qubit requirements in this scenario. Regarding the number of jobs, equation \ref{eq:qubitrequire} shows that the number of qubits increases linearly with the number of jobs. To conclude, our algorithm requires a reasonable amount of qubits (between one hundred and several thousand for real-size industrial instances).

\subsubsection{Circuit depth}

The circuit depth analysis is more complex due to the potential future capabilities of machines to natively implement advanced gates and the inherent challenge of optimizing quantum circuit depth, which can lead to high constants when circuits are not optimized. The depth of the fixed-point quantum search algorithm is primarily determined by the depth of the oracles, which is multiplied by the number of iterations of the algorithm.

As mentioned in section \ref{subsubsec:ResourceConstraintOracle}, the resource constraint oracle has a theoretical depth of $\mathcal{O}\left[\frac{I(I-1)}{2}K\log(CK)\right]$. If we use the full search approach, the depth of the time window and spacing constraint oracle is dominated by $\mathcal{O}(\log((C-1)K+1))^2$, as stated in section \ref{subsubsec:FeasiblePathOracle}. Meanwhile, the quantum walk-inspired scheme has a depth dominated by $\mathcal{O}(K \log(CK)^2)$ (section \ref{subsubsec:QW}). Overall, the circuit depth will asymptotically increase quadratically with the number of machines $I$ and follow a logarithmic-linear function with the number of jobs $K$. \\

However, these asymptotic behaviors can obscure potential high-value constants in the circuit depth. Hence, we used the transpile function from the Qiskit library to compute the circuit depth associated with our previous experiments. This transpile function expresses our circuit in terms of the basic single-qubit gates ('u1', 'u2', 'u3') and the C-NOT two-qubit gate ('cx'). For more details, refer to the Qiskit documentation \cite{qiskit2024quantum}. Table \ref{tab:depth} shows the output of the transpile function, with values around forty thousand per iteration for 2 machines and 2 jobs, and twenty million per iteration when we add one machine or one job. These values far exceed the numbers given in the next 3-5 years in the public roadmaps of quantum computing companies.

\begin{table}
    \centering
\begin{tabular}{|c|c|c|}
        \hline
        I & K & Depth per iteration \\
        \hline
        2 & 2 & $4,1$ . $10^5$ \\
        \hline
        3 & 2 & $2,1$ . $10^7$\\
        \hline
        2 & 3 & $2,1$ . $10^7$\\
        \hline
\end{tabular}
\caption{Approximate depth per iteration displayed by the transpile function of the qiskit library to the basic gates ('u1', 'u2', 'u3','cx')}
\label{tab:depth}
\end{table}

These circuit depths almost certainly aren't feasible without employing quantum error correction (QEC). QEC allows exponentially small gate errors (virtually no errors) at the cost of an overhead in the number of qubits, meaning that if we want to use a few hundred error-corrected qubits, the quantum computer will probably have to have a few tens of thousands of physical qubits. Quantum computers of that size do appear in the roadmaps of some quantum computing companies (in the 5-10 year range) \cite{gambetta2024ibm, google2024google, iqm2024development}. While progress is hard to predict, recent results show rapid progress in QEC \cite{paetznick2024demonstration, acharya2024quantum, eberhardt2024pruning}.

\section{Conclusion}

This work aimed to explore the practical use of quantum search in the near future when fault-tolerant quantum machines (FTQC) become available. It originated from the observation that the quadratic speed-up promised by quantum search might not be sufficient to make this approach practical, as the ratio between the number of possible solutions and the number of good solutions typically increases exponentially with the problem size in optimization problems. However, quantum search algorithms do not impose any conditions on the initial state, which means this approach could be relevant if, by exploiting problem structures, we achieve to build an initial superposition for which the aforementioned ratio increases at most quadratically with the problem size.
The idea of state space reduction (SSR) is to exploit the specific constraints structure of an optimization problem to build an initial state superposition of reduced size while on the same time maintaining quantum search for constraints and objective functions where no structure can be exploited. We developed a proof of concept of a state-space reduction algorithm inspired by quantum walk to generate a set of feasible solutions for a specific scheduling problem. Our analysis of the search space size with and without the application of SSR, along with our numerical results on quantum emulators highlight the potential of this promising approach, which could lead to more efficient quantum search processes by focusing on a smaller search space exploiting the problem's structure whenever possible. Our analysis of space size and numerical results demonstrates the potential of this approach: with a similar number of qubits, the initial superposition size grows only quadratically compared to the number of solutions, allowing for fewer iterations, while decreasing the circuit depth and potentially scaling the method to generate solutions for real-size industrial instances in the near future. Specifically, our implementation on a simplified use case using the Qiskit emulator indeed showcased the high potential of the state-space reduction but also pointing out the challenge of the depth requirement for this kind of approach. Future work could explore the practical applications of our technique in maintenance planning problems. In such scenarios, generating a pool of candidate solutions could serve as input for more precise classical algorithms or expert analysis. Additionally, comparing the SSR-Quantum search to classical approaches for this task would be valuable. This method could also be generalized to other problems by extending the quantum walk-inspired scheme to different types of graph structures. These structures often present more complex constraints that challenge classical solutions but can be leveraged to reduce the search space. For instance, this approach could be applied to grid optimization problems, where the grid inherently imposes a structure, or to electric vehicle problems, where positions or charging intervals might be represented within graph structures \cite{veshchezerova2022quantum}. Both of these extensions could lead to significant breakthroughs in solving industrial problems.

\section*{Aknowledgement}
This work was supported by the European project NEASQC (funded from the
European Union’s Horizon 2020 research and innovation programme grant agreement No
951821).






\bibliography{references}
\bibliographystyle{plain}
\end{document}